\documentclass[12pt,twoside,doublespace]{article}

\usepackage{latexsym,amssymb,amsfonts}
\usepackage{graphicx}

\def\smskip{\par\vskip 5 pt}
\def\QED{\hfill $\Box$\smskip}
\newtheorem{theorem}{Theorem}

\newtheorem{proposition}{Proposition}

\newtheorem{definition}{Definition}

\begin{document}

\begin{center}

\vspace{35pt}

{\Large \bf Variational Inequality Type Formulations}

\vspace{5pt}

{\Large \bf  of General Market Equilibrium Problems}

\vspace{5pt}

{\Large \bf  with Local Information}

\vspace{5pt}

\vspace{35pt}

{\sc Igor V.~Konnov\footnote{\normalsize E-mail: konn-igor@ya.ru}}

\vspace{35pt}

{\em  Department of System Analysis
and Information Technologies, \\ Kazan Federal University, ul.
Kremlevskaya, 18, Kazan 420008, Russia.}

\end{center}

\vspace{25pt}

\hrule

\vspace{10pt}

{\bf Abstract:} We suggest a new approach to creation of general market equilibrium
models involving economic agents with local and partial knowledge about the system
and under different restrictions. The market equilibrium problem is then
formulated as a quasi-variational inequality that enables us to establish
existence results for the model in different settings.
We also describe dynamic processes, which fall into
information exchange schemes of the proposed market model.
In particular, we propose an iterative solution method for quasi-variational inequalities,
which is based on evaluations of the proper market information only
in a neighborhood of the current market state without knowledge of the
whole feasible set and prove its convergence.


\vspace{10pt}

\hrule

\vspace{10pt}

{\em MS Classification:} {91B50, 91B55, 49J40, 65K10}



\section{Introduction}

Variational inequalities give a suitable common format for many
applied problems arising in Economics, Mathematical Physics,
Transportation, Communication Systems, Engineering and
 are closely related with other general problems in
Nonlinear Analysis, such as fixed point, game equilibrium, and optimization problems.
For this reason, their theory and solution methods are developed rather well; see e.g.
\cite{BC84,Nag93,Aub98,FP03,Kon13a} and
references therein. Nevertheless, many models describing rather complex behavior
of systems with interdependent elements require a somewhat more general format,
where the fixed feasible set is replaced with
 a set-valued mapping. This modification leads to the so-called quasi-variational inequalities;
see  e.g. \cite{BL84,CP82,BC84,Har91,Aub98,Yua99} and
references therein. However, these problems are much more difficult for solution than the usual
variational inequalities and the existing general solution methods are very expensive and applicable mostly
for small dimensional problems; see e.g. \cite{CP82,HC91,FP05}.

For this reason, we are interested in revealing new general classes of applications, which
can be formulated as (quasi-)variational inequalities in order to derive existence properties
in a unified manner and enable us to develop efficient methods for finding their solution
under rather natural assumptions that may be also treated as dynamic processes converging to equilibrium states
and justifying these equilibrium concepts. In this work, we intend to present
implementation of this approach for general market equilibrium models.
Traditionally, these models play one of the main roles in Economics.
They in particular indicate ways of equilibrating
different interests and opportunities of  active elements (agents,
participants) involved in economic systems and may serve as a basis for
description of behavior of these very complex systems.

We recall that investigation of the general equilibrium models
dates back to the book by L.~Walras \cite{Wal74}.
Since then, a number of different kinds of these models were proposed;
see \cite{Wal33,AD54,SH73,SS75,Aub98}. More detailed surveys and expositions
of the basic contributions in this field can be found in \cite{Nik68,AH71}.
These models usually describe markets of a great number of economic agents (customers and producers)
so that actions of any separate agent can not impact the state of the whole
system, hence any agent utilizes some integral system
parameters (say, prices), rather than
the information about the behavior of other separate agents.
The traditional approach is based on the assumption that the agents
are able to determine precisely their desirable collections of commodities
for any vector of common market prices. Then an equilibrium market state
is defined from the complementarity conditions between the price and excess supply
of each commodity since the usual material balance condition needs additional
restrictive assumptions. At the same time, it holds for each commodity with positive
equilibrium price; see e.g. \cite{Nik68,AH71}.

In this paper, we suggest some other approach to creation of general equilibrium models.
The current market state is supposed to be
defined by the vector of transaction quantities of the agents.
These quantities satisfy the balance equation.
We suppose that the agents have only
partial and local knowledge about the current and future
behavior of this so complicated system and
evaluate their feasible volume transaction and price sets
in a neighborhood of the current market state. Due to these limitations
the general market equilibrium problem
is formulated as a quasi-variational inequality and corresponds to the local
maximal market profit state concept. This approach enables us to utilize well
developed tools from the theory of quasi-variational inequalities for establishing existence
results.
However, finding market equilibrium points by an iterative solution method
seems rather difficult in the general case. We suggest iterative processes,
which have rather natural treatment and are convergent for some classes of
market equilibrium problems. In particular, we propose an iterative solution
method for quasi-variational inequalities,
which is based on evaluations of the market information
in a neighborhood of the current market state rather
than whole feasible set, which may be unknown,
and prove its convergence under rather weak assumptions.

We outline now briefly the further organization of the paper.
In Section \ref{sc:2}, we recall some auxiliary properties and facts from
the theory of quasi-variational inequalities.
 In Section \ref{sc:3}, we describe a single commodity market equilibrium model
 and its variational inequality re-formulation.
In Section \ref{sc:4}, we describe a general multi-commodity market equilibrium model
in the form of a quasi-variational inequality and  give its basic properties.
Implementation issues of this model are discussed in
Section \ref{sc:5}. Section \ref{sc:6} involves description and substantiation of
some dynamic market processes converging to equilibrium points.
In Section \ref{sc:7}, we give a comparison of the presented model
with the existing basic general and partial equilibrium models.
Section \ref{sc:8} contains some conclusions.


\section{Basic Preliminaries}\label{sc:2}

We first recall several continuity properties of set-valued
mappings.


\begin{definition} \label{def:2.1}
{\em   Let $W$ and $V$ be convex sets in the $l$-dimensional Euclidean space $\mathbb{R}^{l}$,  $W \subseteq V$,
 and let $T: V \to \Pi (\mathbb{R}^{l})$ be a set-valued mapping.
 The mapping $T$ is said to be

 (a) {\em upper semicontinuous (u.s.c.)} on $W$,
 if for each point $v \in W$ and for each open set $U$ such
that $U \supseteq T(v)$, there is an open neighborhood $\tilde V $ of $v$ such
that $T(w) \subset U$ whenever $w \in \tilde V  \bigcap W$;

 (b) {\em lower semicontinuous (l.s.c.)} on $W$,
 if for each point $v \in W$ and for each open set $U$ such
that $U \bigcap T(v) \neq \varnothing$, there is an open neighborhood $\tilde V  $ of $v$ such
that $U \bigcap T(w) \neq \varnothing$ whenever $w \in \tilde V  \bigcap W$;

 (c)  {\em continuous} on $W$, if it is both u.s.c. and l.s.c. on $W$.
}
\end{definition}

Here and below, $\Pi (X)$ denotes the power set of $X$,
i.e., the family of all nonempty subsets of $X$.
Let $\tilde H$ be a set in $\mathbb{R}^{l}$
 and let $H: \tilde H \to \Pi (\tilde H)$  and $G: \tilde H \to \Pi (\mathbb{R}^{l})$ be set-valued mappings.
Then we can define the {\em quasi-variational inequality problem} (QVI for short):
Find $x^{*} \in H(x^{*})$ such that
\begin{equation} \label{eq:2.1}
\exists g^{*} \in G(x^{*}), \ \langle g^{*} , y-x^{*} \rangle \geq 0  \quad \forall y \in H (x^{*}).
\end{equation}
If the feasible mapping $H$ is constant, problem (\ref{eq:2.1}) reduces to the usual
{\em variational inequality problem} (VI for short).

It is well known that QVIs can be converted into the fixed point format
by using the projection mapping; see \cite[Theorem 5.1]{CP82}.
Let $\pi _{X} (x)$ denote the projection of  a point $x$ onto a set $X$.

\begin{proposition} \label{pro:2.1}
Suppose that the sets $H(x)$ and $G(x)$ are nonempty,
convex and compact for each $x \in \tilde H$.
Then QVI  (\ref{eq:2.1}) is equivalent to the fixed point
problem
$$
x^{*} \in \pi _{H (x^{*})} [x^{*} -\theta G (x^{*})],
$$
for some  $\theta >0$.
\end{proposition}

We now give existence results for general QVIs on compact sets; see e.g.
\cite[Theorem 9.14]{Aub98}.

\begin{proposition} \label{pro:2.2}
Suppose that the set $\tilde H$ is
convex and compact, the mapping $H : \tilde H \rightarrow \Pi (\tilde H)$
is continuous on  $\tilde H$, the mapping $G : \tilde H \rightarrow \Pi (\mathbb{R}^{l})$
is upper semi-continuous on  $\tilde H$,
the sets $H(x) $ and $G(x)$ are nonempty,
convex and compact for each $x \in \tilde H$.
Then QVI  (\ref{eq:2.1}) has a solution.
\end{proposition}

Moreover, it was noticed in \cite[Theorem 3]{Har91} that the above conditions may be
relaxed if the feasible mapping value includes the current point.

\begin{proposition} \label{pro:2.3}
Suppose that the set $\tilde H$ is nonempty, convex and compact,
the mapping $G : \tilde H \rightarrow \Pi (\mathbb{R}^{l})$
is upper semi-continuous on  $\tilde H$,
the set $G(x)$ is nonempty,
convex and compact for each $x \in \tilde H$,
$x \in H(x)$ for all $x \in \tilde H$.
Then QVI  (\ref{eq:2.1}) has a solution.
\end{proposition}


\section{A single commodity market equilibrium model}\label{sc:3}

For the simplicity of exposition, we begin our considerations
from a simple equilibrium market model of a
homogeneous commodity, which was suggested in
\cite{Kon06,Kon07,Kon07a} for description of
 a simple market based on a suitable auction implementation mechanism.
 Its further development and applications are described in
\cite{Kon09g,Kon13,AGK17}.

The model involves a finite number of traders and buyers of this commodity, their
index sets will be denoted by $I$ and $J$, respectively.
Hence, the roles of agents are fixed.
For each $i \in I$, the $i$-th trader
chooses some offer value $x_{i}$ in his/her capacity segment $[\alpha'_{i}, \beta' _{i}]$ and has a price function $g_{i}$.
Similarly, for each $j \in J$, the $j$-th buyer chooses some bid
value $y_{j}$ in his/her capacity segment $[\alpha'' _{j}, \beta ''_{j}]$
and has a price function $h_{j}$. The signs of all the lower and upper bounds are arbitrary
in general, but the standard choice is to set $\alpha'_{i}=0$ and  $\alpha'' _{j}=0$, hence
the upper bounds are chosen to be positive. Then we can
define the feasible set of offer/bid volumes
\begin{equation} \label{eq:3.1}
 D=\left\{ (x, y) \ \vrule \
\begin{array}{l}
\sum \limits_{i \in I} x_{i} = \sum \limits_{j \in J} y_{j}; \\
x_{i} \in [\alpha' _{i}, \beta' _{i}], i \in I, y_{j} \in [\alpha'' _{j}, \beta ''_{j}], j \in J
\end{array}
\right\},
\end{equation}
where $x = (x_{i})_{i \in I}, y = (y_{j})_{j \in J}$. We suppose
that the prices may in principle depend on offer/bid volumes
of all the participants, i.e. $g_{i}=g_{i}(x, y)$ and $h_{j}=h_{j}(x,
y)$. We say that a pair $(\bar x, \bar y) \in D$ constitutes a
{\em market equilibrium point} if there exists a number $\bar \lambda$ such that
\begin{equation} \label{eq:3.2}
\displaystyle g_{i}(\bar x, \bar y) \left\{
\begin{array}{ll}
\geq
\bar \lambda & \quad \mbox{if} \quad \bar x_{i}=\alpha '_{i}, \\
=\bar \lambda & \quad \mbox{if} \quad \bar x_{i} \in (\alpha
'_{i},\beta'_{i}), \\
\leq \bar \lambda & \quad \mbox{if} \quad \bar x_{i}=\beta'_{i},
\end{array}
\right. \quad \mbox{for} \ i \in I,
\end{equation}
and
\begin{equation} \label{eq:3.3}
\displaystyle h_{j}(\bar x, \bar y) \left\{
\begin{array}{ll}
\leq \bar \lambda & \quad \mbox{if} \quad \bar y_{j}=\alpha''_{j}, \\
=\bar \lambda & \quad \mbox{if} \quad \bar y_{j} \in (\alpha''_{j},
     \beta ''_{j}), \\
\geq \bar \lambda & \quad \mbox{if} \quad \bar y_{j}=\beta ''_{j},
\end{array}  \right. \quad \mbox{for} \ j \in J.
\end{equation}
Observe that the number $\bar \lambda$ can be treated as
a market clearing price, which equilibrates the market and
yields also the offer/bid volumes for all the participants.
In fact, the minimal offer (bid)
volumes correspond to traders (buyers) whose prices are greater
(less) than $\bar \lambda$, and the maximal offer (bid) volumes correspond
to traders (buyers) whose prices are less (greater) than $\bar \lambda$.
The prices of other participants are equal to $\bar \lambda$ and their
volumes may be arbitrary within their capacity bounds, but should be
subordinated to the balance equation. It follows that agents' prices at equilibrium
may in general differ from the market clearing price.

In \cite{Kon06} (see also \cite{Kon07,Kon07a}), the following basic
relation between the equilibrium problem
(\ref{eq:3.1})--(\ref{eq:3.3}) and a VI was established.

\begin{proposition} \label{pro:3.1} \hfill \\
(a) If $(\bar x, \bar y, \bar \lambda)$ satisfies
(\ref{eq:3.2})--(\ref{eq:3.3}) and $(\bar x, \bar y) \in D$, then
$(\bar x, \bar y)$ solves VI: Find $(\bar x, \bar y) \in D$ such
that
\begin{equation}\label{eq:3.4}
\sum \limits_{i \in I} g_{i} (\bar x, \bar y) (x_{i} - \bar x_{i}) -
\sum \limits_{j \in J} h_{j} (\bar x, \bar y) (y_{j} - \bar y_{j})
\geq 0 \quad \forall (x, y) \in D.
\end{equation}
(b) If a pair $(\bar x, \bar y)  \in D$ solves VI (\ref{eq:3.4}), then
there exists $\bar \lambda$ such that $(\bar x, \bar y, \bar \lambda)$
satisfies (\ref{eq:3.2})--(\ref{eq:3.3}).
\end{proposition}
Therefore, we can apply various results from the theory of VIs (see e.g. \cite{Kon07})
for investigation and solution of the above equilibrium problem.
For instance, if the set $D$ in (\ref{eq:2.1})
is nonempty and bounded,  the functions $g_{i}$ and $h_{j}$ are continuous
for all $i \in I$ and $j \in J$, then VI (\ref{eq:3.4}) has a solution.


\section{A general multi-commodity market equilibrium model} \label{sc:4}

We now present a multi-commodity extension of the model
described in Section \ref{sc:3}, thus also extending those in
\cite{Kon07a,Kon13}.

The model is an $n$-commodity market involving $m$ economic
agents, they are producers and consumers with respect to these
commodities.  We suppose that the agents can strike preliminary (or virtual) bargains
until an equilibrium state will be attained. Then all these deals
are fixed. Denote by $N=\{1, \dots, n\}$ and
$I=\{1, \dots, m\}$ the index sets of commodities and agents at this market.
Observe that roles of the agents with respect to the commodities at the market may be changed.
Next, let the vector $x^{i}=(x_{i1}, \ldots, x_{in})^{\top}$ define the current
transaction quantities of commodities of the $i$-th agent, so that $x_{ij}>0$ means
that his/her current  sold volume of the $j$-th commodity equals $x_{ij}$, whereas
$x_{ij}<0$ means that his/her current  purchased volume of the $j$-th commodity
equals $-x_{ij}$. The current market state is thus supposed to be completely
defined by the virtual volume vector $x=(x^{i})_{i \in I} $.

Behavior of economic agents may in general depend on their current various
goals and restrictions. In our model,
the agents determine first their current feasible volume
transaction sets $Y_{i}(x) \subseteq \tilde Y_{i}$ for
$i \in I$, which are supposed to be nonempty sets in $\mathbb{R}^{n}$ and attributed to each market state $x$.
Similarly, the agents determine  their feasible price sets $P_{i}(x) \subseteq \tilde P_{i}$ for
$i \in I$, which are also supposed to be nonempty sets in $\mathbb{R}^{n}_{+}$ and attributed to each market state $x$.
The sets $\tilde Y_{i}$ and $\tilde P_{i}$, $i \in I$ give the total market transaction bounds
(respectively, the total market price bounds) that may be unknown to the agents in general.
These sets define the set-valued market transaction mapping
$x \mapsto Y(x) $ on the set $\tilde Y$ where
$$
Y(x)=\prod_{i \in I} Y_{i}(x), \quad \tilde Y=\prod_{i \in I} \tilde Y_{i},
$$
and the set-valued feasible market price mapping
$x \mapsto P(x) $ on the set $\tilde Y$ where
$$
P(x)=\prod_{i \in I} P_{i}(x), \quad P(x) \subseteq \tilde P=\prod_{i \in I} \tilde P_{i} \quad \forall x \in \tilde Y.
$$
Both volume and price mappings describe the behavior of the economic
agents with respect to market states.
Using these sets, the agents can make some transactions and change their current transaction volumes.
Next, any market transaction must also satisfy the balance equation, hence we
obtain the set-valued feasible market transaction mapping $x \mapsto D(x) $ with the values
\begin{equation}\label{eq:4.1a}
   D (x) =\left\{ y \in \mathbb{R}^{nm} \ \vrule \
\sum \limits_{i \in I} y^{i} = \mathbf{0};  \ y^{i} \in Y_{i}(x), \ i \in I
\right\}.
\end{equation}
Taking into account Proposition \ref{pro:3.1}
we say that a vector $\bar x \in D (\bar x)$ is a
{\em market equilibrium point} if
\begin{equation} \label{eq:4.2a}
\exists \bar p \in P (\bar x), \ \langle \bar p , y-\bar x \rangle \geq 0  \quad \forall y \in D (\bar x).
\end{equation}
Therefore, the general multi-commodity market equilibrium problem
is formulated as a QVI. We observe that the value
$$
-\langle \bar p, \bar x \rangle =- \sum \limits_{i \in I}\langle \bar p^{i} ,\bar x^{i} \rangle
$$
gives precisely the difference between the total bought and sold amount at
the current market state $\bar x$ and current price vector $\bar p$ of the agents,
i.e., it is the current possible profit of the market. Therefore,
the market equilibrium point $\bar x \in D (\bar x)$ can be treated as the
state that provides the maximal possible market profit at the
current prices of the agents subject to the
constraints determining the balance and current feasible transactions.

We now give analogues of
equilibrium conditions (\ref{eq:3.2})--(\ref{eq:3.3}) for the above problem.
We will take the following basic assumptions.

{\bf (A1)} {\em The sets $\tilde Y_{i} \subset \mathbb{R}^{n}$ are
convex and compact for all $i \in I$. At each state $x \in \tilde Y$
the sets $Y_{i}(x)$ and $P_{i}(x)$ are nonempty,
convex and compact for all $i \in I$.}

Let us define the problem of finding a feasible point $\bar x \in D (\bar x)$ and a
vector $\bar \lambda \in \mathbb{R}^{n}$ such that
\begin{equation} \label{eq:4.1}
\exists \bar p \in P (\bar x), \ \sum \limits_{i \in I}\langle \bar p^{i} -\bar \lambda,y^{i}-\bar x^{i} \rangle  \geq 0
  \quad \forall y \in Y (\bar x).
\end{equation}

\begin{proposition} \label{pro:4.1} Suppose that the assumptions  in {\bf (A1)} are fulfilled. Then:

(a)  If a point $\bar x  \in D (\bar x)$ solves QVI (\ref{eq:4.1a})--(\ref{eq:4.2a}), then
there exists $\bar \lambda$ such that $(\bar x, \bar \lambda)$
satisfies (\ref{eq:4.1}).

(b) If a pair $(\bar x, \bar \lambda) \in D (\bar x) \times \mathbb{R}^{n}$ satisfies
(\ref{eq:4.1}), then
$\bar x$ solves QVI (\ref{eq:4.1a})--(\ref{eq:4.2a}).
\end{proposition}
{\bf Proof.}
If $\bar x  \in D (\bar x)$ solves QVI (\ref{eq:4.1a})--(\ref{eq:4.2a}), then it
solves the optimization problem
\begin{equation}\label{eq:4.2}
{\min} \ \to \ \left\{ \langle \bar p, y \rangle \ | \  y \in D (\bar x) \right\}.
\end{equation}
By using the suitable Karush-Kuhn-Tucker theorem for this problem
(see e.g. \cite[Corollary 28.2.2]{Roc70}), we obtain that
(\ref{eq:4.1}) gives its necessary and sufficient optimality condition. Conversely,  let
 a pair $(\bar x, \bar \lambda) \in D (\bar x) \times \mathbb{R}^{n}$ satisfy
(\ref{eq:4.1}).  Using the same Karush-Kuhn-Tucker theorem
we obtain that $\bar x $ solves (\ref{eq:4.2}), hence it is a solution to
QVI (\ref{eq:4.1a})--(\ref{eq:4.2a}). \QED

We observe that problem (\ref{eq:4.1}) can be replaced by the following equivalent system of
partial QVIs:
\begin{equation} \label{eq:4.3}
\exists \bar p^{i} \in P _{i}(\bar x), \ \langle \bar p^{i} -\bar \lambda,y^{i}-\bar x^{i} \rangle
\geq 0  \quad \forall y^{i} \in Y _{i}(\bar x), \quad \forall i \in I.
\end{equation}

It is clear that (\ref{eq:4.1}) (or (\ref{eq:4.3})) represents an analogue of
the equilibrium conditions in (\ref{eq:3.2})--(\ref{eq:3.3}) for the above
multi-commodity market equilibrium model, hence the
point $\bar \lambda$ in (\ref{eq:4.1}), which is nothing
but the Lagrange multiplier vector corresponding to the
balance constraint, should be treated as the genuine
market equilibrium  price vector for this model. We observe that
the equilibrium prices of the agents do not in general coincide with
$\bar \lambda$ due to the presence of the constraints
determining  the feasible volume transaction sets $Y_{i}(\bar x)$
for all $i=1,\ldots,l$. This is the case even for the simplest single commodity
equilibrium model due to (\ref{eq:3.2})--(\ref{eq:3.3}).
It follows that the agents can not maintain in general the unique
market prices at a given state of the market as in the usual Walrasian models.

We  give the existence result for
the market equilibrium problem as proper specialization of
similar properties for general QVIs from Section \ref{sc:2}.

{\bf (A2)} {\em The mapping $Y : \tilde Y \rightarrow \Pi (\mathbb{R}^{nm})$
is continuous on  $\tilde Y$, the mapping $P : \tilde Y \rightarrow \Pi (\mathbb{R}^{nm})$
is upper semi-continuous on  $\tilde Y$.}


\begin{theorem} \label{thm:4.1}
Suppose that assumptions  {\bf (A1)} and {\bf (A2)} are fulfilled.
Then QVI  (\ref{eq:4.1a})--(\ref{eq:4.2a}) has a solution.
\end{theorem}

The assertion follows from Proposition \ref{pro:2.2}.
Moreover, it seems now natural to suppose that
the feasible volume transaction set of each $i$-th agent $Y_{i}(x)$
will contain the current volume vector $x^{i}$.
This assumption allows us to take the existence result for the usual
VIs as  in Proposition \ref{pro:2.3}. Set
\begin{equation}\label{eq:4.4a}
  \tilde D =\left\{ y \in \tilde Y  \ \vrule \
\sum \limits_{i \in I} y^{i} = \mathbf{0} \right\};
\end{equation}
 cf. (\ref{eq:3.1}).

{\bf (A2${}'$)} {\em The mapping $P : \tilde Y \rightarrow \Pi (\mathbb{R}^{nm})$
is upper semi-continuous on  $\tilde Y$, $x \in Y(x)$ for all $x \in \tilde D$, $\tilde D \neq \varnothing$.}


\begin{theorem} \label{thm:4.2}
Suppose that assumptions  {\bf (A1)} and {\bf (A2${}'$)} are fulfilled.
Then QVI  (\ref{eq:4.1a})--(\ref{eq:4.2a}) has a solution.
\end{theorem}
{\bf Proof.} Let us take the following VI: Find
a point $\bar x \in \tilde D$  such that
$$
\exists \bar p \in P (\bar x), \ \langle \bar p , y-\bar x \rangle \geq 0  \quad \forall y \in \tilde D.
$$
Under the above assumptions it must have a solution; see e.g.
\cite[Theorem 9.9]{Aub98}. Due to {\bf (A2${}'$)} it
also solves  QVI (\ref{eq:4.1a})--(\ref{eq:4.2a}). \QED


\section{Implementation issues} \label{sc:5}

The multi-commodity equilibrium model presented in Section \ref{sc:4} seems
incomplete without a more detailed description of
feasible transaction and price sets attributed to
each state $x$ of the market. We suppose that the agents have only
partial and local knowledge about the  market system.
Hence, they are able to only evaluate their
market, industry and social goals and restrictions
in a neighborhood of the current market state.
That is, each agent may evaluate his/her temporal capacity bounds by using the information
about the current supply or demand of some related commodities as well as
their bounds. For this reason, the feasible transaction set
$Y_{i}(x)$ of each $i$-th economic agent is determined by his/her
market capacity constraints.
These constraints can reflect both personal and common market or technology
relationships among several commodities that can be actively
utilized by different agents.  Also, the initial distribution of
commodity endowments must have clear impact on the feasible transaction sets.
For instance, the sets $Y_{i}(x)$ can be determined as intersections
of personal and common obligatory feasible sets.
Let $\tilde Y_{i} \subset \mathbb{R}^{n}$ denote the
set of obligatory personal capacity constraints
of the $i$-th agent, which is independent of market states.
Also, let $W_{i} \subset \mathbb{R}^{nm}$
denote the set of obligatory common market volume constraints
of the $i$-th agent. Then the feasible transaction set
$Y_{i}(x)$ at the state $x$ is defined as follows:
$$
    Y_{i}(x)=\tilde Y_{i} \bigcap Y'_{i}(x), \ \mbox{where} \
    Y'_{i}(x)=\left\{ y^{i} \in \mathbb{R}^{n} \ \vrule \ (x^{-i},y^{i}) \in W_{i} \right\}.
$$
In the general case, the sets $W_{i}$ may also depend on the current state $x$, i.e.
their precise values may be unknown.
Similarly, each agent may insert also the set of desirable market volume distributions
in the above definition. But in the simplest case
$Y_{i}(x)$ is a box-constrained set of the form
\begin{equation} \label{eq:5.1}
Y_{i}(x)=\left\{ y^{i} \in \mathbb{R}^{n} \ \vrule \ y_{ij} \in [\alpha _{ij}(x), \beta _{ij}(x)], \ j =1,\ldots,n \right\},
\end{equation}
where $\alpha _{ij}(x)$ and $ \beta _{ij}(x)$ give the current estimates of
market capacity bounds by the $i$-th economic agent.

The evaluation of the current feasible price set $P_{i}(x)$ by the $i$-th economic agent
can be based on possible industrial capacity and social constraints
and the agent's profit as the goal function.
In other words, each $i$-th agent solves the optimization problem:
\begin{equation}\label{eq:5.2}
\max \limits _{p^{i} \in V_{i}(x)} \to  \langle  p^{i} , x^{i} \rangle,
\end{equation}
where $V_{i}(x) \subset \mathbb{R}^{n}_{+}$ and takes its solution set as $P_{i}(x)$.
We note that $\langle  p^{i} , x^{i} \rangle$
gives precisely the profit of the $i$-th agent within his/her own prices.
The feasible set $V_{i}(x)$ in (\ref{eq:5.2}) may depend on the current state $x$
and reflect also utilization of common waste treatment plants
or environment restrictions for some industry technologies of different agents,
partial budget type restrictions, common transportation capacity restrictions,
application of nonlinear production technologies, utilization of special financial tools, etc.

Let us take the simplest case with independent linear production technologies.
For each $x^{i}$ define the index sets $I'_{i}$ and $I''_{i}$,
which determine the agent's supply and demand commodities.
It follows that $x_{ij} \geq 0$ if $j \in I'_{i}$ and
$x_{ij} \leq 0$ if $j \in I''_{i}$.
Then $V_{i}(x) \equiv V_{i}$ becomes a polyhedral set, i.e.
\begin{equation} \label{eq:5.3}
V_{i}=\left\{ v^{i} \in \mathbb{S}^{n}_{+} \ \vrule \  \sum _{s\in I''_{i}} a^{i} _{sj} v_{is} \geq v_{ij}, \ j\in I'_{i} \right\},
\end{equation}
where $a^{i} _{sj}$ is the amount of the $s$-th commodity used for production of one unit of the $j$-th commodity by the
$i$-th economic agent,
$$
\mathbb{S}^{n}_{+}=\left\{ z \in \mathbb{R}^{n}_{+} \ \vrule \  \sum _{j =1}^{n} z _{j}=1 \right\}
$$
is the standard simplex in $\mathbb{R}^{n}$. Hence, the relative price scale is taken here.
Then (\ref{eq:5.2})--(\ref{eq:5.3}) reduces to the linear programming
problem:
\begin{equation}\label{eq:5.3a}
\max \limits _{p^{i} \in V_{i}} \to  \langle  p^{i} , x^{i} \rangle,
\end{equation}
besides, we have
\begin{equation} \label{eq:5.4}
P_{i}(x)=\partial \mu _{i}(x^{i}),
\end{equation}
i.e. the solution set of problem (\ref{eq:5.3})--(\ref{eq:5.3a}) is the sub-differential
of its optimal value function $\mu _{i}(x^{i})=\mu _{i}(x)$, which is clearly convex.
The inequalities in (\ref{eq:5.3}) mean that the minimal price $p_{ij}$
of the $j$-th commodity can not be greater than the expenses for all the factors (resources)
used for production of one unit of this commodity.
We observe that the non-industrial consumption commodities can be taken into account similarly.
The agent should only indicate  that consumption expenses per unit of the $s$-th commodity
will be covered with selling all the commodities such that $a^{i} _{sj}>0$.
The coefficient values $a^{i} _{sj}$ will give the corresponding relative scale, i.e.
the relative weights of commodities for any agent.
In such a way this model can in principle describe a pure exchange market.
Similarly, common waste treatment or environment restrictions can be inserted in the model.
In the general case all these parameters and the index sets $I'_{i}$ and $I''_{i}$
may also depend on the current state $x$, i.e.
the agents can change their preferences for commodities with respect to the current
distribution of volumes.


\section{Convergence of dynamic processes} \label{sc:6}

The existence of natural dynamic market processes converging to
an equilibrium point is very essential for justification
of any suggested equilibrium concept. But it is well known that
substantiation of such dynamic processes is much more difficult in comparison with
that of static existence results. In fact, the multi-commodity equilibrium model
from Section \ref{sc:4} is formulated as QVI  (\ref{eq:4.1a})--(\ref{eq:4.2a}).
However, even finding its solution by an iterative solution method
seems very difficult in the general case; see e.g. \cite{CP82,HC91,FP05}.
In this case, we have additional restrictions due to
partial and local knowledge of the agents about the system.
In fact, we intend to describe a process of changing the market states
that continues until attaining a solution of problem (\ref{eq:4.1a})--(\ref{eq:4.2a}), thus creating a
{\em decentralized (self-regulation) transaction mechanism}.
The implementation of this mechanism is clearly attributed to
a suitable information exchange scheme of this model and depends on
the properties of the defined sets and mappings.
Therefore, any chosen iterative solution method must admit
a natural treatment of this method as a
dynamic process for finding market equilibrium points in the above model
together with suitable convergence conditions.
For this reason, it is more suitable to select rather broad classes of
equilibrium problems of form (\ref{eq:4.1a})--(\ref{eq:4.2a}) and suitable
natural converging dynamic processes. First of all we observe that
converging dynamic processes based on bilateral transactions were proposed
for the single commodity market equilibrium model from Section \ref{sc:3};
see \cite{Kon15e,Kon16b,Kon19}. We now give examples of iterative processes,
which have rather natural treatment and are convergent for some classes of
multi-commodity market equilibrium problems.

Let us take assumptions {\bf (A1)}  and {\bf (A2${}'$)} and suppose in addition that
the set $P_{i}(x)$ is determined as the solution set
of optimization problem (\ref{eq:5.3})--(\ref{eq:5.3a})
for all $i \in I$ and $x \in \tilde Y$.
Then, due to (\ref{eq:5.4}), QVI (\ref{eq:4.1a})--(\ref{eq:4.2a}) becomes equivalent
to the quasi-optimization problem: Find a point $\bar x \in D (\bar x)$ such that
\begin{equation} \label{eq:6.1a}
\mu(\bar x) \geq \mu(y) \quad \forall y \in D (\bar x),
\end{equation}
where
$$
\mu(x)=\sum \limits_{i \in I} \mu _{i}(x^{i}),
$$
and $\mu _{i}(x^{i})=\mu _{i}(x)$ is the optimal value of problem (\ref{eq:5.3})--(\ref{eq:5.3a}).
This means that the set of market equilibrium states coincide with the set of the points
with minimal pure expenses (or maximal profit) of the whole market.
The appearance of the moving feasible set $D(x)$ in (\ref{eq:6.1a}) reflect both relationships
and restricted knowledge of the agents.
Hence, problem (\ref{eq:6.1a}) can be also treated as a relative optimization problem; see \cite{Kon19d}.

Bearing in mind Proposition \ref{pro:2.1}, we can  apply the simplest projection method
\begin{equation} \label{eq:6.2a}
x^{(k+1)} = \pi _{D (x^{(k)})} [x^{(k)}-\theta_{k} p^{(k)}], \ p^{(k)} \in P(x^{(k)}), \ \theta_{k}>0,
\ k=0,1,\ldots,
\end{equation}
for finding a solution of problem (\ref{eq:6.1a}).
However, its convergence needs additional assumptions since
problems (\ref{eq:5.3})--(\ref{eq:5.3a}) may have many solutions, i.e.
the function $\mu $ is in general non-differentiable. This gives a quasi-optimization
problem with convex non-differentiable function and creates certain difficulties for
convergence properties of method (\ref{eq:6.2a}). At the same time, extensions of more sophisticated
non-differentiable optimization methods (see e.g. \cite{Sho79,DV81,Kon13a})
do not admit a natural treatment within the above
market equilibrium model and its information exchange scheme.
We can take the stationary case where the agents utilize only fixed feasible transaction sets.
More precisely, we take the following assumptions.

{\bf (B1)} {\em At each state $x \in \tilde Y$ it holds that $Y_{i}(x)=\tilde Y_{i}$, where
the set $\tilde Y_{i} \subset \mathbb{R}^{n}$ is convex and compact for each $i \in I$. The set
$\tilde D$ in (\ref{eq:4.4a}) is non-empty.}

{\bf (B2)} {\em At each state $x \in \tilde Y$ the set $P_{i}(x)$ is determined as the solution set
of optimization problem (\ref{eq:5.3})--(\ref{eq:5.3a}), where the set $V_{i}$ is non-empty
for all $i \in I$.}

Then QVI (\ref{eq:4.1a})--(\ref{eq:4.2a}) reduces to the set-valued VI:
Find a vector $\bar x \in D$ such that
\begin{equation} \label{eq:6.2}
\exists \bar p \in P (\bar x), \ \langle \bar p , y-\bar x \rangle \geq 0  \quad \forall y \in \tilde D.
\end{equation}
Moreover, due to (\ref{eq:5.4}) this problem is equivalent to the convex optimization problem:
$$
\min \limits _{x \in \tilde D} \to  \mu(x);
$$
cf. (\ref{eq:6.1a}).
It is clear that VI  (\ref{eq:6.2}), (\ref{eq:4.4a}) has a solution under the assumptions
in {\bf (B1)}--{\bf (B2)}. This result follows e.g. from Theorem \ref{thm:4.1}.
Next, we can find a solution of VI  (\ref{eq:6.2}), (\ref{eq:4.4a})  by using the fixed point
iterate (\ref{eq:6.2a}) that now reduces to the sub-gradient projection method:
\begin{equation} \label{eq:6.4}
x^{(k+1)} = \pi _{\tilde D} [x^{(k)}-\theta_{k} p^{(k)}], \ p^{(k)} \in \partial \mu(x^{(k)}), \ \theta_{k}>0,
\ k=0,1,\ldots,
\end{equation}
with the standard step-size rules:
\begin{equation} \label{eq:6.5}
\sum \limits_{k=0}^{\infty }\theta _{k}=\infty , \ \sum
\limits_{k=0}^{\infty }\theta ^{2}_{k}<\infty .
\end{equation}
Its convergence property is well known and can be deduced e.g. from Lemma 2.1 in  \cite[Chapter V]{GT89}.


\begin{proposition} \label{pro:6.1}
Suppose that  assumptions  {\bf (B1)}--{\bf (B2)} are fulfilled. If a sequence $\{x^{(k)}\}$ is subordinated to
rules (\ref{eq:6.4})--(\ref{eq:6.5}), it converges to a solution of problem (\ref{eq:6.2}), (\ref{eq:4.4a}).
\end{proposition}

Nevertheless, implementation of each iterate may be not so easy within
the information exchange scheme of the multi-commodity market model.
Bearing in mind (\ref{eq:5.1}), we specialize the market capacity bounds as follows.

{\bf (B1${}'$)} {\em At each state $x \in \tilde Y$
it holds that $Y_{i}(x)=\tilde Y_{i}$, where
$$
\tilde Y_{i}=\left\{ y^{i} \in \mathbb{R}^{n} \ \vrule \ y_{ij} \in [\alpha _{ij}, \beta _{ij}], \ j =1,\ldots,n \right\}
$$
for all $i \in I$. The set
$\tilde D$ in (\ref{eq:4.4a}) is non-empty.}

Clearly, {\bf (B1${}'$)}--{\bf (B2)} imply {\bf (B1)}--{\bf (B2)}. According to
the sub-gradient projection method, the agents first determine
their price vector $p^{(k)}$ for the current market state $x^{(k)}$
from the optimization problems of form (\ref{eq:5.3})--(\ref{eq:5.3a}). Afterwards, they find the next state $x^{(k+1)}$
from (\ref{eq:6.4}). This problem now decomposes into $n$ independent
single commodity market equilibrium problems of form (\ref{eq:3.4}) with
affine price functions. In fact, the $j$-th commodity market problem
is written as VI: Find $x^{(k+1)}_{(j)}\in \tilde D_{j}$ such that
\begin{equation}\label{eq:6.6}
\sum \limits_{i \in I} (p^{(k)}_{ij} +\theta ^{-1}_{k}
(x^{(k+1)}_{ij}-x^{(k)}_{ij})) (y_{ij} - x^{(k+1)}_{ij})
\geq 0 \quad \forall y_{(j)} \in \tilde D_{j},
\end{equation}
where $y_{(j)}=(y_{1j}, \ldots, x_{mj})^{\top}$  and
\begin{equation} \label{eq:6.7}
\tilde D_{j}=\left\{ y_{(j)} \in \mathbb{R}^{m} \ \vrule \
\sum \limits_{i \in I} y_{ij}= 0, \ y_{ij} \in [\alpha _{ij}, \beta _{ij}], \ i \in I \right\}.
\end{equation}
This problem can be solved by simple solution algorithms including the bilateral exchanges;
see \cite{Kon15e,Kon16b}. Therefore, process (\ref{eq:6.4}) can be naturally implemented
within usual market mechanisms.

We still intend to present the simple method
for finding a solution of QVI (\ref{eq:4.1a})--(\ref{eq:4.2a})
within the assumptions from {\bf (A1)}  and {\bf (A2${}'$)} with moving
feasible transaction sets.

{\bf (C1)} {\em The sets $\tilde Y_{i} \subset \mathbb{R}^{n}$ are
convex and compact for all $i \in I$. At each state $x \in \tilde Y$
the set $Y_{i}(x) \subseteq \tilde Y_{i}$ is  convex and compact for each $i \in I$.}

In general, the set-valuedness of the price mappings  $P_{i}(x)$ may
create certain difficulties for the agents since the clear choice of unique prices is more suitable.
For this reason, we suppose that
each $i$-th economic agent bearing in mind some
most suitable reference price vector $\bar v^{i}$ chooses his/her
feasible price value $P_{i}(x)$ at state $x$, i.e. problem (\ref{eq:5.3})--(\ref{eq:5.3a})
is replaced by the optimization problem
\begin{equation}\label{eq:6.8}
\max \limits _{p^{i} \in V_{i}} \to  \sum \limits^{n}_{j=1}
\left\{ p_{ij}x_{ij} -0.5 \beta_{i} (p_{ij} - \bar v_{ij})^{2}\right\},
\end{equation}
where $\beta_{i}>0$ is a suitable weight parameter chosen by the $i$-th agent.
Since it has the unique solution $p^{i} (x^{i} )=p^{i} (x)$,
the mapping $x \mapsto P_{i}(x) $  is now single-valued for each $i\in I$ and so is $x \mapsto P(x) $.
If we denote by $\eta _{i}(x^{i} )$ the optimal value of problem (\ref{eq:6.8}), (\ref{eq:5.3}), then
the function $\eta _{i}(x^{i} )$ will be convex and differentiable and $\eta' _{i}(x^{i} )=p^{i} (x^{i} )$.
In turn, $p(x)=(p^{i} (x))_{i\in I}$ is the gradient of the convex function
\begin{equation}\label{eq:6.10}
\eta(x)=\sum \limits_{i \in I} \eta _{i}(x^{i}).
\end{equation}

In other words, we take the following assumptions.

{\bf (C2)} {\em At each state $x \in \tilde Y$ the point $p^{i}(x)=p^{i} (x^{i} )$ is determined as a unique solution
of optimization problem (\ref{eq:6.8}), (\ref{eq:5.3}), where the set $V_{i}$ is non-empty
for all $i \in I$.}

{\bf (C3)} {\em At each state $x \in \tilde D$ it holds that $x \in Y(x)$,
the set $ \tilde D$  is non-empty,
the mapping $D : \tilde Y \rightarrow \Pi (\tilde D)$ is lower semi-continuous  on  $\tilde Y$.}

Then QVI (\ref{eq:4.1a})--(\ref{eq:4.2a}) reduces to the single-valued problem:
Find a vector $\bar x \in D (\bar x)$ such that
\begin{equation} \label{eq:6.11}
\langle p(\bar x) , y-\bar x \rangle \geq 0  \quad \forall y \in D (\bar x).
\end{equation}
Its solution can be however deduced from the usual VI:
Find a vector $\bar x \in \tilde D$ such that
\begin{equation} \label{eq:6.12}
\langle p(\bar x) , y-\bar x \rangle \geq 0  \quad \forall y \in \tilde D,
\end{equation}
where the set $\tilde D$ is defined in (\ref{eq:4.4a}).
From {\bf (C1)} and {\bf (C3)} it follows that $\tilde D$ is nonempty,
convex and compact. From {\bf (C2)} it follows that the mapping $p(x) $ is continuous
 on $\tilde Y$. Hence, VI (\ref{eq:6.12}) has a solution $\bar x \in \tilde D$.
Due to {\bf (C3)} we have  $\bar x \in Y(\bar x)$, therefore,  $\bar x \in D (\bar x)$.
This means that $\bar x$ is a solution to QVI (\ref{eq:6.11}). We observe that the reverse implication
(\ref{eq:6.11}) $\Longrightarrow $ (\ref{eq:6.12}) does not hold in general.


\begin{proposition} \label{pro:6.2}
Suppose that  assumptions  {\bf (C1)}--{\bf (C3)} are fulfilled.
Then QVI  (\ref{eq:6.11}) has a solution.
\end{proposition}

Moreover, VI (\ref{eq:6.12}) is equivalent to the convex smooth optimization problem:
$$
\min \limits _{x \in \tilde D} \to  \eta(x),
$$
where $\eta(x)$ is defined in (\ref{eq:6.10}).
It follows that we can in principle apply one of various well known
smooth optimization methods (see e.g. \cite{PD75}) for
finding a solution of QVI  (\ref{eq:6.11}).
However, most of them are too sophisticated and do not fall into the above
market information exchange framework. For instance, even simple projection
conditional gradient methods with line-search or
utilization of a priori information can not be taken due to the
indicated information restrictions. For this reason, we will present a dynamic process based
on the parametric conditional gradient method without line-search,
which was proposed in \cite{Kon18d} for custom optimization problems.
Given a point $x \in \tilde D$ we define the auxiliary problem
\begin{equation} \label{eq:6.14}
 \min \limits _{y \in D(x)} \to \langle  p(x), y\rangle .
\end{equation}
We denote by $Z(x)$ the solution set of problem  (\ref{eq:6.14}),
thus defining the set-valued  mapping $x\mapsto Z(x)$. Observe that
$Z(x)$ is determined only by the local information from $D(x)$ rather than $\tilde D$
and that the set $Z(x)$ is always non-empty, convex, and compact.

\medskip
\noindent {\bf Method (PCGM).} \\
 {\em Initialization:} Choose a point $w^{(0)} \in \tilde D$, a number $\beta \in (0,1)$,
and a positive sequence $\{\delta _{s}\} \to 0$. Set $s=1$.

{\em Step 0:} For the given number $s$ choose a positive sequence $\{\tau _{l,s}\}$
such that $\tau _{l,s} \in (0, 1)$ and $\{\tau _{l,s}\} \to 0$ as $l \to \infty$.
Set  $k=0$, $l=0$, $x^{(0)}=w^{(s-1)}$, and choose a number $\theta_{0} \in (0, \tau _{0,s}]$.

{\em Step 1:} Find a point $y^{(k)}\in Z(x^{(k)})$. If
\begin{equation} \label{eq:6.15}
\langle p(x^{(k)}), x^{(k)}-y^{(k)} \rangle  \geq \delta _{s},
\end{equation}
go to Step 2. Otherwise set $w^{(s)}=x^{(k)}$, $s=s+1$ and go to Step 0. {\em (Restart)}

{\em Step 2:}  Set $d^{(k)}= y^{(k)}-x^{(k)}$, $x^{(k+1)}=x^{(k)}+\theta_{k}d^{(k)}$. If
\begin{equation} \label{eq:6.16}
  \langle p(x^{(k+1)}), d^{(k)}\rangle \leq \beta \langle p(x^{(k)}), d^{(k)}\rangle,
\end{equation}
take $\theta_{k+1} \in [\theta_{k}, \tau _{l,s}]$. Otherwise
take $\theta_{k+1} \in (0,\min \{\theta_{k}, \tau _{l+1,s}\}]$ and  set $l=l+1$.
Afterwards set $k=k+1$ and go to Step 1.
\medskip

Note that each outer iteration (stage) in $s$ contains some number of inner iterations in $k$
with the fixed tolerance $\delta _{s}$. Completing each stage,
that is marked as restart, leads to decrease of its value.
Note that the choice of the parameters $\{\tau _{l,s}\}$ can be in principle independent
for each stage $s$. Also, by (\ref{eq:6.15}), we  have
$$
\langle p(x^{(k)}),d^{(k)}\rangle \leq -\delta_{s}<0
$$
in (\ref{eq:6.16}). It follows that
 \begin{eqnarray}
 \eta (x^{(k+1)})-\eta  (x^{(k)}) &\leq & \theta_{k} \langle p(x^{(k+1)}),d^{(k)}\rangle
 \leq \beta \theta_{k} \langle p(x^{(k)}),d^{(k)}\rangle \nonumber \\
 &\leq &  -\beta \theta_{k}\delta_{s}, \label{eq:6.17}
\end{eqnarray}
besides, $x^{(k)} \in \tilde D$ for each $k$.

We show that each stage is well defined.


\begin{proposition} \label{pro:6.3}
Suppose that  assumptions  {\bf (C1)}--{\bf (C3)} are fulfilled.
Then the number of iterations at each stage $s$ is finite.
\end{proposition}
{\bf Proof.}
Fix any $s$ and suppose that the sequence $\{x^{(k)}\}$ is infinite.
Then the number of changes of index $l$ is also infinite.
In fact, otherwise we have $\theta _{k} \geq \bar \theta>0$ for $k$ large enough, hence
(\ref{eq:6.17}) gives
$$
 \eta (x^{(k+t)}) \leq \eta (x^{(k)})-t \beta \bar \theta \delta_{s} \ \to -\infty \ \mbox{as} \  t \to \infty,
$$
for $k$ large enough, which is a contradiction.
Therefore, there exists an infinite subsequence of indices $\{k_{l}\}$
 such that
\begin{equation} \label{eq:6.18}
\langle p(x^{(k_{l}+1)}), d^{(k_{l})}\rangle > \beta \langle p(x^{(k_{l})}), d^{(k_{l})}\rangle,
\end{equation}
where $d^{(k_{l})}=y^{(k_{l})}-x^{(k_{l})}$. Besides, it holds that
$$
 \theta_{k_{l}} \in (0,\tau_{l,s}], \ \theta_{k_{l}+1} \in (0,\tau_{l+1,s}],
$$
where
$$
\lim \limits_{l\rightarrow \infty}\tau_{l,s}=0.
$$
Both the sequences $\{x^{(k)}\}$ and $\{y^{(k)}\}$ belong to the bounded
set $\tilde D$ and hence have limit points. Without loss of generality, we can
suppose that the subsequence $\{x^{(k_{l})}\}$ converges to a point
$\bar x$ and the corresponding subsequence $\{y^{(k_{l})}\}$ converges to a point
$\bar y$. Due to (\ref{eq:6.15}) we have
\begin{equation} \label{eq:6.19}
\langle p(\bar x),\bar y -\bar x \rangle= \lim \limits_{l\rightarrow \infty}
\langle p(x^{(k_{l})}), y^{(k_{l})}-x^{(k_{l})} \rangle \leq -\delta_{s}.
\end{equation}
At the same time, taking the limit $l\rightarrow \infty$ in (\ref{eq:6.18}), we obtain
$$
\langle p(\bar x),\bar y -\bar x \rangle \geq  \beta \langle p(\bar x),\bar y -\bar x
\rangle,
$$
i.e.,  $   (1-\beta ) \langle p(\bar x),\bar y -\bar x \rangle \geq 0$,
which is a contradiction with (\ref{eq:6.19}). \QED
We are ready to prove convergence of the whole method.


\begin{theorem} \label{thm:6.1}
Suppose that  assumptions  {\bf (C1)}--{\bf (C3)} are fulfilled.
Then:

(i)  The number of changes of index $k$ at each stage $s$  is finite.

(ii) The sequence $\{w^{(s)}\}$ generated by method (PCGM) has limit points, all
these limit points are solutions of QVI  (\ref{eq:6.11}).

\end{theorem}
{\bf Proof.} Assertion (i) has been obtained in Proposition \ref{pro:6.3}.
By construction, the sequence $\{w^{(s)}\}$ is bounded, hence it has limit points.
By definition, for each $s$ it holds that
\begin{equation} \label{eq:6.20}
\langle p(w^{(s)}),y-w^{(s)}\rangle \geq -\delta_{s} \quad \forall y \in D (w^{(s)}).
\end{equation}
Take an arbitrary limit point $\bar w$ of $\{w^{(s)}\}$, then
$$
\bar w =\lim \limits_{t\rightarrow \infty }w^{(s_{t})},
$$
for some subsequence $\{w^{(s_{t})}\}$, hence $\bar w \in \tilde D$. Then
$D(\bar w )$ is nonempty and $\bar w \in D(\bar w)$. Take any
$\bar y \in D(\bar w)$, then there exists a sequence of points $\{
y^{(s_{t})} \} \to \bar y$, $ y^{(s_{t})}\in D(w^{(s_{t})})$,   since the mapping $D$ is lower semi-continuous on  $\tilde Y$.
Setting $s=s_{t}$ and $y=y^{(s_{t})}$ in  (\ref{eq:6.20}) and taking the limit $t\rightarrow \infty$,
we obtain
$$
\langle p(\bar w),  \bar y-\bar w \rangle \geq 0.
$$
This means that $\bar w$ is a solution of QVI  (\ref{eq:6.11}).
Assertion (ii) is true.
\QED

Again, implementation of each iterate may be not so easy in general.
Bearing in mind (\ref{eq:5.1}), we specialize the feasible transaction sets
and market capacity bounds as follows.

{\bf (C1${}'$)} {\em At each state $x\in \tilde Y$
it holds that
$$
Y_{i}(x)=\left\{ y^{i} \in \tilde Y_{i} \ \vrule \ y_{ij} \in [\alpha _{ij}(x), \beta _{ij}(x)], \ j =1,\ldots,n \right\},
$$
where
$$
\tilde Y_{i}=\left\{ y^{i} \in \mathbb{R}^{n} \ \vrule \ y_{ij} \in [\alpha' _{ij}, \beta' _{ij}], \ j =1,\ldots,n \right\}
$$
for all $i \in I$, the set
$$
\tilde Y=\prod_{i \in I} \tilde Y_{i}
$$
is bounded.}

Clearly, {\bf (C1${}'$)}, {\bf (C2)}--{\bf (C3)} imply {\bf (C1)}--{\bf (C3)}.
Following (PCGM), the agents determine
their price vector $p^{(k)}$ for the current market state $x^{(k)}$
in accordance with {\bf (C2)}. However, changing the stage does not require
calculation of the price since the real state remains the same.
This means that the agents only adjust their current tolerances for determination
of the sufficient decrease of the pure expenses of the market; see (\ref{eq:6.17}).
The main part of the iteration consists in solution of problem  (\ref{eq:6.14})
at $x^{(k)}$.  This problem clearly decomposes into $n$ independent
single commodity market equilibrium problems of form (\ref{eq:3.4}) with
fixed prices. In fact, the $j$-th commodity market problem is now written as follows:
$$
\min \limits _{y_{(j)} \in D_{j}} \to  \sum \limits_{i \in I} p^{(k)}_{ij} y_{ij},
$$
where
$$
D_{j}=\left\{ y_{(j)} \in \mathbb{R}^{m} \ \vrule \
\sum \limits_{i \in I} y_{ij}= 0, \ y_{ij} \in [\alpha^{(k)} _{ij}, \beta^{(k)} _{ij}], \ i \in I \right\},
$$
$\alpha^{(k)} _{ij}$ and $ \beta^{(k)} _{ij}$ are the corresponding capacity bounds; cf. (\ref{eq:6.6})--(\ref{eq:6.7}).
This problem can be solved by simple arrangement algorithms. Therefore,
(PCGM) can be naturally implemented within usual market mechanisms.

It should be also noticed that (PCGM) has rather low computation expenses per
iteration and provides opportunities for computations in a distributed manner where
 computational units possess only limit data sets and information
 transmission flows are also limited. These multi-agent computational
 procedures are developed very intensively; see e.g. \cite{BT89,SFSPP14}.


\section{Relationships with some other basic equilibrium models}\label{sc:7}

In this section, we intend to make a comparison of the presented model
with the existing basic general and partial equilibrium models.

We recall that the classical Walrasian type general equilibrium models describe a
market of a great number of economic agents so that actions of any separate agent
can not impact the state of the whole system and traditionally belong to the perfect competition market models.
These models are based on the assumption that
any agent does not utilize particular information about the behavior of the others,
but accepts the common market price values and then determines precisely his/her personal
demand and/or supply values. These values create the general market excess supply values
that have certain impact on the prices. Namely, a negative (positive) market excess supply of a commodity
forces its price to increase (decrease), thus defining the so-called t\^{a}tonnement
process; see \cite{Wal74,AH71}. Hence, the model essentially exploits (in fact,
postulates) the assumption that there exists a common (market) price value for each commodity at any moment
that can be recognized by any separate agent.
Then the budget constraint of any consumer involves just these common prices,
which is usually active for his/her optimal choice due to the basic non-satiability
property of the utility function. Similarly, any producer maximizes his/her profit,
which is again determined by the same common prices.
Next, the equilibrium conditions in the classical models are determined as
complementarity relationships between the prices and
market excess supply values; see e.g. \cite{Nik68,AH71}.

Our model is based on the assertion that the agents of a general market
do not possess in fact a sufficient information about the current and future
behavior of this so complicated system. That is, information for their decision making can be only
partial and local. This means that each agent can evaluate his/her opportunities only
in a neighborhood of the current market state by using his/her technology peculiarities and capacities,
some integral parameters of the market (say, the current total supply or demand for certain commodities),
and some information on other agents due to the transaction process.
For instance, this lack of information
about the market participants is typical for contemporary
telecommunication based systems; see e.g. \cite{IK10,RAR10}.
Besides, the behavior of each economic agent depends on
his/her current market, industry and social goals and restrictions.
For this reason, the assumption that there exists an \lq\lq ideal broker"
who is able to determine precisely the common market price for each commodity at any moment
and report simultaneously all these prices to all the agents seems rather restrictive.
Similarly, it seems too difficult requirement for an agent to precisely fulfil the budget constraint
with common prices also because of the data uncertainty.
It should be noted that common constraints may appear due to the complexity of the system,
these constraints can reflect common market, transportation, and/or technology
relationships that are actively utilized by different agents.
For instance, markets with joint constraints arise often in
telecommunication and energy sectors; see e.g. \cite{HH04,IK10,PSPF10}.

For this reason, we suppose that the agents will utilize their own prices that may differ from the current
market prices even at an equilibrium state. These prices are values of their feasible price mappings
attributed to the current market state that is determined by the current distribution of commodities.
Instead of the exact budget constraint any agent should simply indicate current sources that will cover his/her
consumption expenses for each purchased commodity, as described in Section \ref{sc:5}.
Similarly, the agents determine feasible transaction sets dependent of the current market state.
This approach leads to a quasi-variational inequality formulation of the market equilibrium state, which
corresponds to the local maximal market profit in a neighborhood of the current state.
It was proved in  Section \ref{sc:4} that the proposed market equilibrium model also involves
common market prices as the Lagrange multipliers corresponding to the common
balance constraint, but they are in general different from the equilibrium prices of the agents;
see Proposition \ref{pro:4.1}. It also follows that the agents may have positive endowments at equilibrium
due to the market capacity bounds as in the single commodity case.

We observe that most imperfectly competitive (or partial) equilibrium models (see e.g. \cite{OS90})
and some general market equilibrium models (see e.g. \cite{SS75,Yua99})
are formulated as non-cooperative game-theoretic problems.
Due to the common balance constraint in any market, the agents are mutually dependent,
which contradicts to the custom non-cooperative game setting with equal and independent players.
But in the imperfect competition case, actions of each agent
can change the state of the whole market so that agents' utility functions depend
on these actions (strategies).
In order to avoid the mentioned drawback, the imperfect competition models are based on
uneven roles of market sides. Namely, one side of the market (say, consumers)
is presented by a common price (inverse demand) function, whereas each agent from the other side (producer)
can change the state of the whole system by choosing his/her supply
so that the common balance is always satisfied implicitly.
Therefore, the model becomes a game of producers and its solution is usually determined as
its Nash equilibrium point. However, this is not the case for the general market equilibrium models
with equal agents, which are then formulated as extended non-cooperative game-theoretic problems;
see e.g. \cite{SS75,Yua99}.  The presence of the joint binding constraints creates certain
difficulties for players and needs in a special non-trivial mechanism for attaining such an equilibrium point;
see e.g. \cite{Kon16c}. Therefore, our quasi-variational inequality formulation of the market equilibrium
problem seems rather natural and suitable for different applications.


\section{Conclusions}\label{sc:8}

We described a new class of general market equilibrium
models involving economic agents with restricted knowledge about the system.
Due to these limitations the market equilibrium problem
was formulated as a quasi-variational inequality. We then deduced existence
results for the model from the theory of quasi-variational inequalities.
Besides, we indicated some ways of implementation of the proposed model
in different settings and made comparisons with the other basic general  equilibrium models.
We also presented decentralized dynamic
processes converging to equilibrium points within
information exchange schemes of the market model.
In particular, we proposed an iterative solution method for quasi-variational inequalities,
which is based on evaluations of the market information
in a neighborhood of the current market state rather than whole feasible set, which may be unknown,
and proved its convergence.

Investigations of these models and methods can be continued in several directions.
For instance, it seems worthwhile to adjust them to different applications,
which involve explicit technology, environment  and transportation capacity restrictions
and develop new suitable dynamic processes.


\end{document}